\documentclass{article}
\usepackage{amsfonts}
\usepackage{amsmath}
\usepackage{color}
\usepackage{amssymb}
\usepackage{epsfig}
\RequirePackage{psfrag}
\RequirePackage{graphicx}
\usepackage{multicol}
\usepackage{float}

\usepackage[colorlinks=true]{hyperref}
\hypersetup{urlcolor=blue, citecolor=blue}
\usepackage[top = 2cm, left=3cm]{geometry}

\newtheorem{theorem}{\bf Theorem}[section]
\newtheorem{lemma}[theorem]{\bf Lemma}
\newtheorem{proposition}[theorem]{\bf Proposition}

\newtheorem{definition}[theorem]{\bf Definition}

\newtheorem{remark}[theorem]{{\bf Remark}.\it}{}
\newenvironment{keywords}{\small {\bf key words}.\it}{\vskip 10pt}
\newenvironment{AMS}{\small {\bf AMS subject classification}.\it}{\vskip 10pt}
\newenvironment{proof}{{\it Proof: }}{\hfill $\square$}

\def\O{\Omega}
\def\eps{\varepsilon}
\def\OiT{\Omega_i\times(0,T)}
\def\Qq{\mathcal{Q}}
\def\Dd{\mathcal{D}}

\def\phii{\varphi_i}
\def\ds{\displaystyle}
\def\R{\mathbb{R}}
 \def\N{\mathbb{N}}
\def\l{\lambda}
 \def\a{\alpha}
 \def\b{\beta}
 
 \def\beqn{\begin{eqnarray}}
 \def\eeqn{\end{eqnarray}}
 \def\be{\begin{equation}}
 \def\ee{\end{equation}}
 \def\nn{\nonumber}
 \def\refe#1{(\ref{#1})}
\def\ue{u^\eps}

\def\iint{\int\!\!\!\int}
 \def\Mm{\mathcal{M}}
 \def\Rr{\mathcal{R}}

 \def\k{\kappa}
 
 \def\t{\tilde}

\def\Cc{\mathcal{C}}

\begin{document}
\author{Cl\'ement  Canc\`es\thanks{UPMC Univ Paris 06, UMR 7598, Laboratoire Jacques-Louis Lions, F-75005, Paris, France (\href{mailto:cances@ann.jussieu.fr}{\tt cances@ann.jussieu.fr})}
\thanks{The author is partially supported by GNR MoMaS}}
\title{Asymptotic behavior of two-phase flows in heterogeneous porous media for capillarity depending only on space.\\ {II}. Non-classical shocks to model oil-trapping}
\date{\today}
\maketitle

\begin{abstract}
We consider a one-dimensional problem modeling two-phase flow in heterogeneous porous media made of two homogeneous subdomains, with discontinuous capillarity at the interface between them. 
We suppose that the capillary forces vanish inside the domains, but not on the interface. 
Under the assumption that the gravity forces and the capillary forces are oriented in opposite directions, 
we show that the limit, for vanishing diffusion, is not in general the optimal entropy solution of the hyperbolic scalar conservation law as in the first paper of the series~\cite{NPCX}. A non-classical shock can occur at the interface, modeling oil-trapping. 
\end{abstract} 

\section{Introduction}

The models of two-phase flows provide good first approximations to predict the motions of oil in the subsoil. Although the theoretical knowledge concerning the question of the existence and the uniqueness of the solution to such models for homogeneous porous media 
\cite{AKM90,CJ86} and for media with regular enough variations \cite{Chen01} is quite complete, few results are available for discontinuous media, as for example media made of several rock types \cite{ABE96, BPvD03, BLS09, FVbarriere, CGP09,  EEM06}. 
\vskip 10pt
One says that oil-trapping occurs when some oil can not pass through interfaces between 
different rocks. Such a phenomenon plays an important role in the basin modeling, to predict the position of eventual reservoirs where oil could be collected. 
As already explained in \cite{BPvD03,vDMdN95}, discontinuities of the capillary pressure field can induce the so-called oil-trapping phenomenon.
\vskip 10pt

The effects of capillarity, which play a crucial role in oil-trapping, seem to play a less important role concerning the motion of oil in homogeneous porous media, and can sometime be neglected to provide the so-called Buckley-Leverett equation.
\vskip 10pt

In this paper, we show that even if the dependence of the capillary pressure with respect to the oil-saturation of the fluid vanishes, the capillary pressure field still plays a crucial role to determine the saturation profile. In order to carry out this study, we restrict our frame to the one-dimensional case.
We will strongly use some recent results \cite{BLS09, FVbarriere,CGP09} obtained on flows in heterogeneous media with discontinuous capillary forces.

\vskip 10pt

We consider a one-dimensional porous medium, made of two different rocks, represented by $\O_1=\R_-^\star$ and $\O_2 = \R_+^\star$.
Let $\pi(u,x)$ be the capillary pressure, then it it is well known (see e.g. the introduction of the associate paper~\cite{NPCX}) 
that, if both phases have different densities, the equation governing the two phase flow can be written
\be\label{dep_NC}
\partial_t u + \partial_x \Big( q c(u,x) + g(u,x) \left(1 - C\partial_x \pi(u,x)\right)  \Big) = 0,
\ee
where $u$ is the \emph{oil saturation} of the fluid, $q$ is the total flow rate, supposed to be a nonnegative constant, $C$ is a constant depending on the buoyancy forces and
$$c(u,x) = c_i(u), \quad g(u,x)=g_i(u), \quad \textrm{ and } \quad
\pi(u,x)=\pi_i(u) \quad \textrm{ if }x\in \O_i.$$
The functions $c_i$ are supposed to be increasing and Lipschitz continuous with $c_i(0)=0$ and $c_i(1) = 1$, while $g_i$ are supposed to be Lipschitz continuous, strictly positive in $(0,1)$ satisfying $g_i(0)=g_i(1) = 0$ and $\pi_i$ are increasing Lipschitz continuous functions.
\vskip 10pt

Physical experiments suggest that the dependence 
of $\pi_i$ with respect to $u$ can be  weak, at least for $u$ far from $0$ and $1$. 
So we want to choose $\pi_1(u)=P_1$, and $\pi_2(u)=P_2$. The equation \refe{dep_NC} 
turns formally to the scalar conservation law with discontinuous flux function
\begin{equation}\label{eq:base}
\partial_t u + \partial_x  f(u,x) =0,
\end{equation}
where $f(u,x)$ (resp. $f_i(u)$) is equal to $q c(u,x) - g(u,x)$ (resp. $q c_i(u) - g_i(u)$).
\vskip 10pt

Such conservation laws have been widely studied in the last years. For a large overview on this topic, we refer to the introduction 
of \cite{BKT09}, or in a lesser extent to the associated paper \cite{NPCX}. In particular, it has been proven by Adimurthi, Mishra and Veerappa Gowda \cite{AMV05} that there might exist an infinite number of $L^1$-contraction semi-groups corresponding to the equation~\eqref{eq:base}. 
Among them, in the case where the functions $f_i$ have at most a single extremum in (0,1), we mention the so-called \emph{optimal entropy solution} which corresponds to the unique entropy solution in the case of a continuous flux function $f_1=f_2=f$. We refer to \cite{AMV05} and to the first part of this communication~\cite{NPCX} for a discussion on the so-called \emph{optimal entropy condition}.

\begin{figure}[htb]
\begin{center}
\resizebox{8cm}{!}{
\input{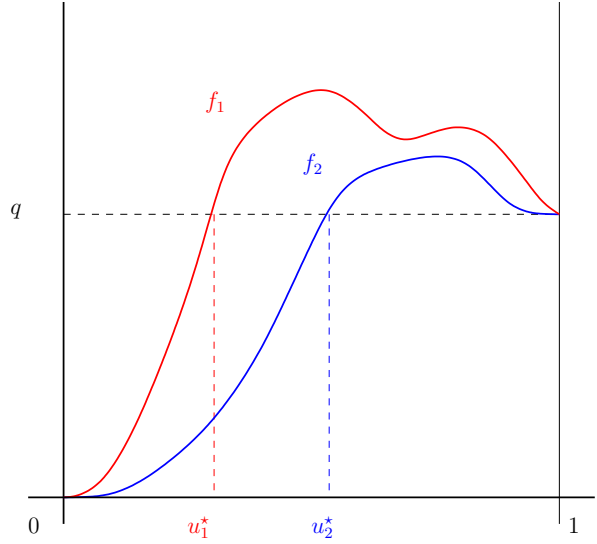}
}
\caption{example of functions $f_i$ satisfying Assumption~{\bf (H1)}. Note the we have note supposed, as it is done in~\cite{AJV04,BKT09}, that $f_i$ has a single local extremum in $(0,1)$, but all the extrema have to be strictly greater than $q$.}\label{fig:f_i}
\end{center}
\end{figure}

\vskip 10pt
In the sequel of this paper, we suppose that
\vskip 10pt
\begin{itemize}
\item[{\bf (H1)}] for $i\in\{1,2\}$,  there exists a value $u_i^\star\in [0,1)$ such that $f_i(u_i^\star) = q$, $f_i$ is increasing on $[0,u_i^\star]$ and $f_i(s) >q $ for all $s\in (u_i,1)$.
\end{itemize}
\vskip 10pt
 We refer to Figure~\ref{fig:f_i} for an illustration of the previous assumption.
We denote by $$\phii(u)= C \int_0^u g_i(s)\, ds.$$
For technical reasons, we have to assume that 
\vskip 10pt
\begin{itemize}
\item[{\bf (H2)}] there exist $R>0$, $\alpha>0$ and $m\in(0,1)$ such that 
\be f_1\circ\varphi_1^{-1}(s)\ge q+ R (\varphi_1(1)-s)^m \quad \textrm{ if } s\in [\varphi_1(1)-\a,\varphi_1(1)].  \label{hyp:Holder}
\ee
\end{itemize}
\vskip 10pt
These assumptions are fulfilled by models widely used by the engineers, 
for which a classical choice of $c_i, g_i$ is 
$$
c_i(u) = \frac{u^{a_i}}{u^{\a_i} + \frac{a}{b} (1-u)^{\b_i}}, \qquad g_i(u)=K_i \frac{u^{\a_i} (1-u)^{\b_i}}{b u^{a_i} + a (1-u)^{\beta_i}},
$$
where $\a_i, \beta_i \ge 1$ and $a,b$ are given constants.
\vskip 10pt

The goal of this paper is to show that if the capillary forces at the level of the interface $\{x=0\}$ are 
oriented in the opposite sense with respect to the gravity forces (in our case $P_1<P_2$), then 
a \emph{non classical} stationary shock can occur at the interface. It was shown by Kaasschieter~\cite{Kaa99} 
that if the capillary pressure field is continuous at the interface (corresponding to the case $P_1 = P_2$), then 
the good notion of solution is the one of \emph{optimal entropy solution}, computed by Adimurthi, Jaffr\'e and Veerappa Gowda 
using a Godunov-type scheme~\cite{AJV04}. 
We have pointed out in \cite{NPCX} that if the capillary forces and the gravity forces are oriented in the same sense, the 
good notion of solution is also the one of \emph{optimal entropy solution}. If the assumptions stated above are fulfilled, if $P_1<P_2$ and if the initial data $u_0$ is large enough to ensure that both phases move in opposite directions, i.e. 
\be\label{hyp:large_init}
u_i^\star \le u_0(x) \le 1 \quad \textrm{ a.e. in }\O_i, 
\ee
 we will show that the limit is not the optimal entropy solution, but the entropy solution to the problem
\be\label{Plim}\tag{$\mathcal{P}_{\lim_{}}$}
\left\{
\begin{array}{l}
\partial_t u + \partial_x f_i(u)=0, \\
u(x=0^-)=1 \textrm{ and } u(x=0^+)=u_2^\star,  \\
u(t=0)=u_0.
\end{array}
\right.
\ee
In the sequel, we denote by $a^+ $ (resp. $a^-$) the positive (resp. negative) 
part of $a$, i.e. $\max(0,a)$ (resp. $\max(0,-a)$), and for $i=1,2$, for $u,\k\in[0,1]$, one denotes 
by 
$$
\Phi_{i+}(u,\k)= \left\{\begin{array}{ll}
f_i(u)-f_i(\k) &\textrm{ if } u\ge \k,\\
0 &\textrm{ otherwise, }
\end{array}\right.
$$
$$
\Phi_{i-}(u,\k)= \left\{\begin{array}{ll}
f_i(\k)-f_i(u) &\textrm{ if } u\le \k,\\
0 &\textrm{ otherwise, }
\end{array}\right.
$$
and 
$$
\Phi_i(u,\k)= \Phi_{i+}(u,\k)+\Phi_{i-}(u,\k)=f_i(\max(u,\k))-f_i(\min(u,\k)).
$$
We can now define the notion of solution to \refe{Plim}, which is in fact an entropy 
in each subdomain $\O_i$, with an internal boundary condition at the level of 
the interface.
\vskip 10pt
\begin{definition}[solution to \refe{Plim}]\label{Sol_NC_def}
Let $u_0\in L^\infty(\R)$, $u_i^\star \le u_0(x) \le 1$ a.e. in $\O_i$, 
A function $u$ is said to be a solution of \refe{Plim} if it belongs to $L^\infty(\R\times\R_+)$, 
$u_i^\star \le u \le 1$ a.e. in $\OiT$, and for $i=1,2$, for all $\psi\in \Dd^+(\overline\O_i\times\R_+)$, for all $\k\in[0,1]$, 
\beqn
\lefteqn{ \int_{\R_+} \int_{\O_i}(u(x,t)-\k)^\pm \partial_t \psi dxdt + \int_{\O_i}(u_0(x)-\k)^\pm \psi(x,0) dx} \nn \\
&&\ds \hspace{10pt}+ \int_{\R_+} \int_{\O_i} \Phi_{i\pm}(u(x,t),\k) \partial_x \psi(x,t) dxdt 
+ M_{f_i} \int_{\R_+} (\overline u_i -\k)^\pm \psi(0,t) dt \ge 0, \label{sol_NC_for}
\eeqn
where $M_{f_i}$ is a Lipschitz constant of $f_i$, and $\overline u_1=1$,  $\overline u_2=u_2^\star$.
\end{definition}
\vskip 10pt
For a given $u_0$ in $L^\infty(\R)$, there exists a unique solution $u$ to \refe{Plim} 
in the sense  of Definition \ref{Sol_NC_def}, which is in fact made on an apposition of two 
entropy solutions in $\R_\pm\times\R_+$. We refer to \cite{MNRR96,Otto96_CRAS} and \cite{Vov02} for proofs of existence and uniqueness to solutions to the problem~\eqref{Plim}. Moreover, thanks to \cite{Cont_L1}, one can suppose that $u$ belongs to $\Cc(\R_+;L^1_{loc}(\R))$.
\vskip 10pt
\begin{theorem}\label{thm:well-posed_Plim}
Let $u_0\in L^\infty(\R)$ with $u_i^\star \le u_0 \le 1$ a.e. in $\O_i$,  then there exists a unique solution to \eqref{Plim} in the sense of Definition~\ref{Sol_NC_def}. Furthermore, if $v$ is another solution to \eqref{Plim} corresponding to $v_0\in L^\infty(\R)$ with $u_i^\star \le v_0 \le 1$ a.e. in $\O_i$, then for all $\R>0$, for all $t\in \R_+$
\be\label{eq:L1_contract_Plim}
\int_{-R}^R (u(x,t) - v(x,t))^\pm dx  \le \int_{-R - M_f t}^{R + M_f t}(u_0(x) - v_0(x))^\pm dx
\ee
where $M_f$ is a Lipschitz constant of both $f_i$.
\end{theorem}
\vskip 10pt
Assume now that both phases move in the same direction:
\be\label{hyp:small_init}
0 \le  u_0(x)  \le u_i^\star \quad \textrm{ a.e. in }Ê \O_i, 
\ee
then it will be shown that the relevant solution $u$ to the problem is the unique entropy solution defined below.
\vskip 10pt
\begin{definition}\label{Def:entro}
A function $u$ is said to be an entropy solution if it belongs to $L^\infty(\R\times\R_+)$, 
$0\le u \le u_i^\star$ a.e. in $\OiT$, and for $i=1,2$, for all $\psi\in \Dd^+(\R\times\R_+)$, for all $\k\in[0,1]$, 
\beqn
\lefteqn{ \int_{\R_+} \int_{\R}\left| u(x,t)-\k \right| \partial_t \psi dxdt + \int_{\R} \left|u_0(x)-\k\right| \psi(x,0) dx} \nn \\
&&\ds + \int_{\R_+}\hspace{-3pt}\sum_{i\in\{1,2\}} \hspace{-3pt}\int_{\O_i}\hspace{-3pt} \Phi_{i}(u(x,t),\k) \partial_x \psi(x,t) dxdt
+  |f_2(\k) - f_1(\k)|\hspace{-3pt}\int_{\R_+}\hspace{-3pt} \psi(0,t) dt \ge 0. \label{eq:entro_for_Tow}
\eeqn
\end{definition}
Thanks to Assumption~{\bf (H1)}, there exist no $\chi\in[0, \max u_i^\star]$ such that $f_1(\chi) = f_2(\chi)$, $f_1$ is decreasing and $f_2$ is increasing on $(\chi-\delta, \chi+\delta)$ for some $\delta>0$. Then the notion of entropy solution described by \eqref{eq:entro_for_Tow} introduced by Towers \cite{Tow00,Tow01} is equivalent to the notion of \emph{optimal entropy solution} introduced in \cite{AMV05} (see also \cite{BKT09}). We take advantage of this by using the very simple algebraic relation~\eqref{eq:entro_for_Tow}. 
\vskip 10pt
It has been proven  that the entropy solution $u$ exists and is unique for general flux functions $f$ \cite[Chapters~4 and 5]{Bachmann_These}. In particular, the following comparison and $L^1$-contraction principle holds.
\vskip 10pt
\begin{theorem}\label{thm:entro}
Let $u_0\in L^\infty(\R)$ with $0 \le u_0 \le u_i^\star$ a.e. in $\O_i$,  then there exists a unique entropy solution in the sense of Definition~\ref{Def:entro}. Furthermore, if $v$ is another entropy solution corresponding to $v_0\in L^\infty(\R)$ with $0 \le v_0 \le u_i^\star$ a.e. in $\O_i$, then for all $\R>0$, for all $t\in \R_+$
\be\label{eq:L1_contract_entro}
\int_{-R}^R (u(x,t) - v(x,t))^\pm dx  \le \int_{-R - M_f t}^{R + M_f t}(u_0(x) - v_0(x))^\pm dx
\ee
where $M_f$ is a Lipschitz constant of both $f_i$.
\end{theorem}

\subsection{non classical shock at the interface}
As already mentioned, the optimal entropy solution can be seen as a extension to the case of discontinuous flux functions of the 
usual entropy solution \cite{K70} obtained for a regular flux function. We will now illustrate that it is not the case with the solution to~\eqref{Plim}. 
Assume for the moment (it will be proved later) that in the case where $u_0(x) \in (u_i^\star,1)$ a.e. in $\O_i$, the corresponding solution $u$ to \eqref{Plim} admits $\overline u_i$ as strong trace on the interface. One has the following \emph{Rankine-Hugoniot} relation
$$
f_1(\overline u_1) = f_2(\overline u_2) = q, 
$$
then $u$ is a weak solution to~\eqref{eq:base}, i.e. it satisfies for all $\psi \in \Dd(\R\times\R_+)$:
\be\label{weak_form}
\int_{\R_+}\int_{\R} u \partial_t \psi\, dxdt+ \int_\R u_0 \psi(\cdot, 0) \, dx + \int_{\R_+}\int_{\R} f(u,\cdot) \partial_x \psi\, dxdt = 0.
\ee

\vskip 10pt
Firstly, suppose for the sake of simplicity that $f_1(u) = f_2(u) = f(u)$, and that $q=0$, then $u_i^\star = 0$ for $i\in\{1,2\}$.
The function 
$$
u(x)  = \left\{ \begin{array}{ll}
1 & \textrm{ if }x<0,\\
0 & \textrm{ if }x>0
\end{array}\right.
$$
is then a steady solution to~\eqref{Plim} satisfying~\eqref{sol_NC_for}. However, since 
$$
\frac{f(1) - f(s)}{1-s} <0 \quad \textrm{ for all } s \in (0,1),
$$
the discontinuity at $\{x=0\}$ does not fulfill the usual \emph{Oleinik entropy condition} (see e.g.~\cite{Smo94}). This discontinuity is thus said to be a 
\emph{non-classical shock}.
\vskip 10pt
Suppose now that $f_1'(1) <0$ and that $f_2'(u_2^\star) >0$, then the pair $(1,u_2^\star)$ is a stationary \emph{undercompressible shock-wave}, that are prohibited for optimal entropy solutions \cite{AMV05} as for classical entropy solutions in the case of regular flux functions. 
\vskip 10pt
\begin{remark}\label{rmk:AB}
It has been pointed out in \cite{AMV05} that allowing a \emph{connection} $(A,B)$, i.e. a stationary undercompressible wave between the left state $A$ and the right state $B$ at the interface lead to another \mbox{$L^1$-contraction} semi-group (see \cite{AMV05,BKT09,GNPT07}), which is so-called entropy solution of type $(A,B)$. However, we rather use the denomination \emph{non-classical shock} for the connection between $A$ and $B$ since, as stressed above, the corresponding solution violates some fundamental properties of the classical entropy solutions. 
\end{remark}
\subsection{oil-trapping modeled by the non-classical shock}\label{subsec:oil-trapping}
In this section, we assume that $q=0$. 
Let $u$ be the solution of the problem \eqref{Plim} corresponding to the initial data $u_0$. Assume that $u$ admits strong traces on the interface. The flow-rate of oil going from $\O_1$ to $\O_2$ through the interface is given by 
$$f_1(\overline u_1) = f_2(\overline u_2) =0.$$ 
Thus the oil cannot overcome the interface from $\O_1$ to $\O_2$, thus if one supposes that $u_0$ belongs to $L^\infty(\R)$, with $0\le u_0 \le 1$ a.e., then the quantity of oil standing between $x=-R$ ($R$ is an arbitrary positive number) and $x=0$ can only grow. 
\vskip 10pt
Indeed, let $t_2>t_1\ge 0$, let $\zeta_n(x) =  \min(1, n(x+R)^+,n x^-)$ and 
$\theta_m(t) = \min(1, m(t-t_1), m(t_2-t))$. Choosing $\psi(x,t)=\zeta_n(x) \theta_m(t)$ in \eqref{weak_form} for $m,n\in \N$ yields, using the positivity of $f_1$
$$
\int_{t_1}^{t_2}\left( \int_{-R}^0u(x,t) \zeta_n(x)dx\right) \partial_t\theta_m(t) dt + \int_{t_1}^{t_2}\theta_m(t) 
\left( \frac1n \int_{-1/n}^0 f_1(u(x,t)) dx \right)dt \le 0.
$$
Since $u$ admits a strong trace on the interface, 
$$
\lim_{n\to\infty} \frac1n \int_{-1/n}^0 f_1(u(x,t)) dx = f_1(\overline u_1) = 0.
$$
Then we obtain 
\begin{equation}\label{ot_1}
\int_{t_1}^{t_2}\left( \int_{-R}^0u(x,t) dx\right) \partial_t\theta_m(t) dt \le 0.
\end{equation}
The solution $u$ belong to $\Cc(\R_+;L^1(\R))$ thanks to \cite{Cont_L1}, thus taking the limit as $m\to \infty$ in \refe{ot_1} provides 
$$
\int_{-R}^0 u(x,t_1)\, dx \le \int_{-R}^0 u(x, t_2)\, dx.
$$
\vskip 10pt
Suppose now that $q\ge0$. Thanks to what follows, we are able to solve the Riemann problem at the interface for any initial data 
$$
u_0(x) = \left\{ \begin{array}{lll}
u_\ell & \textrm{ if }Ê& x<0,\\
u_r & \textrm{ if }Ê& x>0.
\end{array}\right.
$$
The study of the Riemann problem is carried out in Section~\ref{sec:Riemann}, leading to the following result.
\begin{itemize}
\item If $u_\ell> u_1^\star$, then $u_1 = 1$ and $u_2 = u_2^\star$. We obtain the expected non-classical shock at the interface. 
\item If $u_\ell \le u_1^\star$, then $u_1 = u_\ell$ and $u_2$ is the unique value of $[0,u_2^\star]$ such that $f_2(u_2) = f_1(u_\ell)$.
\end{itemize}
Using Assumption~{\bf (H1)}, this particularly implies that in both cases, the flux at the interface is given by 
\be\label{eq:trace_1}
f_1(u_1) = f_2(u_2) = G_1(u_\ell, 1)
\ee
where $G_1$ is the Godunov solver corresponding to the flux function $f_1$:
$$
G_1(a,b) = \left\{\begin{array}{lll}
\ds \min_{s\in[a,b]} f_1(s) & \textrm{ if } & a \le b,\\
\ds \max_{s\in[b,a]} f_1(s) & \textrm{ if } & a > b.
\end{array}\right.
$$
This particularly yields that for any initial data $u_0\in L^\infty(\R)$ with $0\le u_0 \le 1$, the restriction $u_{|_{\O_1}}$ of the solution $u$ to $\O_1$ is the unique entropy solution to 
\be\label{syst:u_1}
\left\{\begin{array}{lll}
\partial_t u + \partial_x f_1(u) = 0 & \textrm{ in } \O_1\times\R_+,\\
u(\cdot, 0) = u_0 & \textrm{ in } \O_1,\\
u(0,\cdot) = \gamma & \textrm{ in }Ê\R_+
\end{array}\right.
\ee
for $\gamma = 1$. Since the solution $u$ to the problem~\eqref{syst:u_1} is a non-decreasing function of the prescribed trace $\gamma$ on $\{x=0\}$, we can claim as in \cite{NPCX} that 
$$
u_{|_{\O_1}} = \sup_{\begin{subarray}{c}\gamma\in L^\infty(\R_+) \\
0 \le \gamma \le 1
\end{subarray}}
\left\{\ v \textrm{ solution to }\eqref{syst:u_1}\ \right\}.
$$
In particular, $u$ is the unique weak solution (i.e. satisfying \eqref{weak_form}) that is entropic in each subdomain and that minimizes the flux through the interface.

\subsection{organization of the paper}
We will introduce a  family of approximate problems  in Section 
\ref{approx_sect_NC}, which takes into account the capillarity, with small 
dependance $\eps$ of the capillary pressure with respect to the 
saturation. We use the transmission conditions introduced in \cite{BLS09,FVbarriere,CGP09,Sch08}  to connect the capillary pressure at the interface.  
For $\eps>0$, the problem~\refe{Pe_NC} admits a unique solution $\ue$ thanks to~\cite{FVbarriere} and it is recalled that a comparison principle holds for the solutions of the approximate problem~\eqref{Pe_NC}. Particular sub- and super-solution are derived in order to show that if $u_0(x) \ge u_I^\star$ a.e. in $\O_i$, then the limit $u$ of the approximate solutions $\left(\ue \right)_{\eps>0}$ as $\eps$ tends to $0$. An energy estimate is also derived.
\vskip 10pt
In Section~\ref{conv_NC_sec}, letting $\eps$ tend to $0$, since no strong pre-compactness can be derived on $\left( \ue \right)_\eps>0$ in $L^1_{loc}(\R\times\R_+)$ from the available estimates, we use the notion of process solution \cite{EGH00}, which is equivalent to the notion of measure valued solution introduced by DiPerna \cite{DP85} (see also \cite{MNRR96,Sze91}). The uniqueness of such a process solution allows us to claim that $(\ue)$ converges strongly in $L^1_{loc}((\R\times\R_+)$ towards the unique solution to \eqref{Plim}. 
\vskip 10pt
In Section~\ref{entro_small}, it is shown that if both phases move in the same direction, that is if $0 \le u_0 \le u_i^{\star}$ a.e. in $\O_i$, then $\left(\ue\right)$ converges towards the unique entropy solution to the problem in the sense of Definition~\ref{Def:entro}.
\vskip 10pt
In Section~\ref{sec:Riemann}, we complete the study of the Riemann problem at the interface.

%
%
\section{The approximate problem}\label{approx_sect_NC}

In this section, we take into account the effects of the capillarity, supposing that 
they are small. We will so build an approximate problem \refe{Pe_NC}, whose unknown 
$\ue$ will depend on a small parameter $\eps$ representing the dependance 
of the capillary pressure with respect to the saturation.
We assume for the sake of simplicity that the capillary pressure in $\O_i$  is 
given by:
 \be\label{pi_NC}
 \pi_i^\eps(\ue)=P_i + \eps \ue.
 \ee
  It has been shown simultaneously in \cite{BLS09} and in \cite{CGP09} that a good way to connect the capillary pressures at the interface is to require 
\be\label{connect_p_NC}
\t \pi_1^\eps(\ue_1)\cap \t \pi_2^\eps(\ue_2) \neq \emptyset,
\ee
where $\ue_1$ and $\ue_2$ are the traces of $\ue$ on the interface, and
 where $\t \pi_i^\eps$ is the monotonous graph given by
$$
\t \pi_i^\eps(s)=\left\{
\begin{array}{ll}
\pi_i^\eps(s) & \textrm{ if } s\in (0,1),\\
(-\infty, P_i] &  \textrm{ if } s=0,\\
{[} P_i+\eps,\infty) &  \textrm{ if } s=1.
\end{array}
\right.
$$

We suppose that the capillary force is oriented in the sense of decreasing $x$, i.e. 
$P_1<P_2$ (the capillary force goes from the high capillary pressure to the low capillary pressure). 
Since $\eps$ is assumed to be a small parameter, we can suppose that 
$0<\eps<P_2-P_1$, so that the relation \refe{connect_p_NC} turns to
\be\label{connect_p_NC2}
\ue_1=1 \textrm{ or }Ê\ue_2=0.
\ee

The flux function in $\O_i$ is then given by:
$$
F_i^\eps(x,t)= f_i(\ue)(x,t)-\eps \partial_x \phii(\ue)(x,t).
$$
Because of the conservation of mass, we require the continuity 
of the flux functions at the interface. Thus the approximate problem becomes
\be\label{Pe_NC}\tag{$\mathcal{P}^\eps$}
\left\{
\begin{array}{l}
\partial_t \ue + \partial_x F_i^\eps=0, \\
\ue(x=0^-)=1 \textrm{ or } \ue(x=0^+)=0,  \\
F_1^\eps (0^-) = F_2^\eps(0^+), \\
u(t=0)=u_0.
\end{array}
\right.
\ee
We are not able to prove the uniqueness of a weak solution of \refe{Plim} if 
the flux $F_i^\eps$ "only" belongs to $L^2(\overline\O_i\times\R_+)$, and we will 
define the notion of prepared initial data, so that the flux belongs to $L^\infty(\O_i\times\R_+)$. 
In this latter case, the uniqueness holds.

\subsection{bounded flux solutions}
We define now the notion of bounded flux solution, that was introduced in this framework in 
\cite{FVbarriere,CGP09}.  
\vskip 5pt
\begin{definition}[bounded flux solution to \refe{Pe_NC}]\label{bounded_Def_NC}
Let $u_0\in L^\infty(\R)$, $0\le u_0 \le 1$, a function $\ue$ 
is said to be a bounded flux solution if 
\begin{enumerate}
\item $\ue \in L^\infty(\R\times\R_+)$, $0\le u \le 1$;
\item $\partial_x \phii(\ue) \in L^\infty(\O_i\times\R_+)\cap L^2_{loc}(\R_+;L^2(\O_i))$;
\item  $u_1^\eps(t)\left(1- u_2^\eps(t)\right) = 0$ for almost all $t\ge 0$, where $u_i^\eps$ denotes the trace of $u^\eps_{|_{\O_i}}$ on $\{x=0\}$.
\item $\forall \psi\in \Dd(\R\times\R_+)$,
\beqn
\lefteqn{ \int_{\R_+}\int_\R \ue(x,t)\partial_t \psi(x,t) dxdt + \int_\R u_0(x)\psi(x,0)dx }\nn \\
&& \ds \hspace{20pt} + \int_{\R_+}\sum_{i\in\{1,2\}}\int_{\O_i} \left[ f_i(\ue) -\eps \partial_x \phii(\ue)    \right]
 \partial_x \psi(x,t) dxdt = 0.    \label{bounded_for_NC}
\eeqn 
\end{enumerate}
\end{definition}
\vskip 10pt
\begin{remark}\label{continuite_temporelle}
Such a bounded-flux $\ue$ solution belongs to $\Cc(\R_+;L^1_{\rm loc}(\R))$, in the sense that there exists 
$\tilde\ue$ in $\Cc(\R_+;L^1_{\rm loc}(\R))$ such that $\ue(t)=\tilde\ue(t)$ for almost all $t\ge 0$ (see \cite{Cont_L1}). More precisely, 
all $t\ge 0$ is a Lebesgue point for $\ue$. So, the slight abuse of notation consisting in considering $\ue(t)$ for all $t\ge 0$ will not lead to any confusion.
\end{remark}
\vskip 10pt
\begin{proposition}\label{comp_prop}
Let $u$ and $v$ be two bounded-flux solutions associated to initial data $u_0,v_0$, 
then for all $\psi\in\Dd^+(\R\times\R_+)$,
\beqn
&\ds \int_{\R_+} \int_\R \left( u -  v \right)^\pm \partial_t \psi dxdt + \int_\R \left( u_0 -  v_0 \right)^\pm \psi(\cdot,0) dx&\nn\\
&\ds + \int_{\R_+} \sum_i \int_{\O_i} \left(  \Phi_{i\pm}(u, v) - \eps \partial_x \left( \phii(u) - \phii(v)\right)^\pm \right)  \partial_x\psi dxdt \ge 0. &\label{comp_prop_eq}
\eeqn
\end{proposition}
\vskip 10pt
We state now a theorem which is a generalization in the case of unbounded domains of Theorem~3.1 and Theorem~4.1 stated in \cite{FVbarriere}. 
\vskip 10pt
\begin{theorem}[existence--uniqueness for bounded flux solutions]\label{comp_bounded_NC}
Let $u_0\in L^\infty(\R)$ with $0 \le u_0 \le 1$ such that:
\begin{itemize}
\item there exists a function $\hat u\in L^\infty(\R)$, with $0 \le \hat u \le 1$ a.e. in $\R$, satisfying $\partial_x \hat u \in L^1(\R)\cap L^\infty(\R)$ and  such that $(u_0-\hat u)\in L^1(\R)$
\item  $\partial_x \phii(u_0)\in L^\infty(\O_i)$ 
\item $\lim_{x\nearrow0} u_0(x) = 1$ or $\lim_{x\searrow0} u_0(x) = 0$.
\end{itemize}
Then there exists a unique bounded flux 
solution $\ue$ to the problem \refe{Pe_NC} in the sense of Definition \ref{bounded_Def_NC} satisfying $(\ue - \hat u) \in L^1(\R)$.
Furthermore, if $\ue, v^\eps$ are two bounded flux solutions associated to initial data
$u_0,v_0$  then 
\be\label{eq:comp_bounded}
u_0(x) \ge v_0(x) \textrm{ a.e. in } \R \quad \Rightarrow \ue(x,t) \ge v^\eps(x,t) \textrm{ a.e. in }\R \textrm{ for all }t\ge 0.
\ee
\end{theorem}
\vskip 10pt
Obviously, the existence of a bounded flux solution can not be extended to any initial data 
in $L^1(\R)$. Indeed, the initial data $u_0$ has at least to involve bounded initial flux, i.e.
$\partial_x \phii(u_0) \in L^\infty(\R)$. An additional natural assumption is needed to 
ensure the existence of such a bounded flux solution : the connection in the graphical 
sense of the capillary pressures at the interface. 
\vskip 10pt
If $(u_0 - \hat u)$ and $(v_0 - \hat u)$ belong to $L^1$ for the same $\hat u$, then \eqref{comp_prop_eq} yields that the bounded flux solutions 
$\ue$ and $v^\eps$ corresponding to $u_0$ and $v_0$ satisfy the following contraction principle: $\forall t\in \R_+$, 
$$
\int_\R (\ue(x,t) - v^\eps(x,t))^\pm dx \le \int_\R (u_0(x) - v_0(x))^\pm dx,
$$
providing the uniqueness result stated in Theorem~\ref{comp_bounded_NC}.

\subsection{particular sub- and super-solutions}\label{sub-super_sec}

We will study particular steady states of the approximate problem~\eqref{Pe_NC}. We will consider 
steady bounded flux solutions $s^\eps$ corresponding to a zero water flow rate, i.e. 
\be\label{eq:steady}
f_i(s^\eps) - \eps \frac{{\rm d} }{{\rm d}x} \phii(s^\eps)= q\quad \textrm{ in } \O_i.
\ee
For $\eps>0$, there are infinitely many solutions $s^\eps$ of the equation~\eqref{eq:steady}. We will 
construct some particular solutions, that will permit us to show that the limit $u$ as $\eps$ tends to $0$  
of bounded flux solutions $\ue$ corresponding to large initial data admits the expected strong traces on 
the interface $\{x=0\}$.
\vskip 10pt

We will introduce now particular solutions of the ordinary differential equation
\be\label{EDO_NC}
y'=f_i\circ\phii^{-1}(y) - q.
\ee
\begin{lemma}\label{EDO_lem}
Let $\eta>0$, there exists a solution $y^\eta$ to \eqref{EDO_NC} for $i=1$ which is nondecreasing on $(-\infty,-0]$ 
equal to $\varphi_1(1)$ on $[-\eta,0)$, satisfying $y^\eta(x) < \varphi_1(1)$ if $x<-\eta$ and $\lim_{x\to-\infty} y^\eta(x) = u_1^\star$.
\end{lemma}
\vskip 10pt
\begin{proof}
Consider the problem 
\be\label{eq:EDO_1}
\left\{\begin{array}{rclll}
w'(x) &=&R(\phii(1)-w(x))^m & \textrm{ if } & x < -\eta, \\
w(-\eta) &=& \varphi_1(1),
\end{array}\right.
\ee
where $R$ and $m$ are constants given by Assumption~{\bf (H2)}. The function 
$$
w^\eta(x) = \phii(1)-\left( R(1-m)(-x-\eta) \right)^{\frac1{1-m}}
$$
is a solution of~\eqref{eq:EDO_1}. Because of {\bf (H2)}, there exists a neighborhood $(-\eta - \delta, -\eta]$ of $\eta$ such that
$w^\eta$ is a super-solution of the problem
\be\label{eq:EDO_2}
\left\{\begin{array}{rclll}
y'(x) &=&f_1\circ\varphi_1^{-1}(y) - q & \textrm{ if } & x < -\eta, \\
y(-\eta) &=& \varphi_1(1).
\end{array}\right.
\ee
Then there exists $y^\eta$ solution to \eqref{eq:EDO_2} such that $y^\eta(x) = \varphi_1(1)$ if $x\in (-\eta,0)$ and 
$$y^\eta(x) \le w^\eta(x) \quad \textrm{ on } (-\eta - \delta, -\eta].
$$
In particular, $y^\eta$ is not constant equal to 1.
Thanks to {\bf (H1)}, the function $y^\eta$ is increasing on the set $\{\ x\in \O_1\ | \ y^\eta(x) \in ( \varphi_1(u_1^\star), \varphi_1(1) )\  \}$.
Assume that there exists $x_\star<-{\eta}$ such that $y^\eta(x_\star) = \varphi_1(u_1^\star)$, then one sets 
$y^\eta(x)  =  \varphi_1(u_1^\star)$ for all $x\in(-\infty, x_\star]$. If $y^\eta(x) > \varphi_1(u_1^\star)$ for all $x<0$, then 
$y^\eta$ is increasing on $(-\infty, -\eta)$. Thus it admits a limit as $x$ tends to $-\infty$, and it is clear that the 
only possible limit is $u_1^\star$.
\end{proof}
\vskip 10pt

\begin{lemma}\label{lem:EDO_2}
Let $\eta>0$, then there exists a solution $z^\eta$ to \eqref{EDO_NC} for $i=2$ which is nondecreasing on $\R$ satisfying 
$z^{\eta}(x) \le \varphi_2\left(\frac{1+u_2^\star}{2}\right)$ if $x\le \eta$, $z^{\eta}(x) \ge \varphi_2\left(\frac{1+u_2^\star}{2}\right)$ if $x\ge \eta$ and 
$\lim_{x\to\infty} z^\eta(x) = \varphi_2(1)$, $\lim_{x\to -\infty} z(x)= u_2^\star$.
\end{lemma}
\vskip 10pt
\begin{proof}
The problem 
$$
\left\{\begin{array}{rclll}
z'(x) &=& f_2\circ\varphi_2^{-1}(z(x)) - q& \textrm{ if } & x\in \R,\\[5pt]
z(\eta) &=&  \varphi_2\left(\frac{1+u_2^\star}{2}\right).
\end{array}\right.
$$
admits a (unique) solution in $\Cc^1(\R,[0,1])$. Since $u_2^\star$ is a constant solution of~\eqref{EDO_NC} for $i=2$, then one has $z(x) \ge u_2^\star$ in 
$\R$. Thanks to {\bf (H1)}, the function $z$ is nondecreasing. This implies that it admits limits respectively in $-\infty$ and in $+\infty$. The only possible values for this limits are respectively $u_2^\star$ and $1$.
\end{proof}
\vskip 10pt
\begin{proposition}\label{prop:sub-super}
Let $\eta>0$, then there exists two families of steady bounded flux solutions $\left(\underline s^{\eps,\eta}\right)_{\eps>0}$ and $\left(\overline s^{\eps,\eta}\right)_{\eps>0}$ tending in $L^1_{loc}(\R)$ as $\eps\to0$ respectively towards 
$$
\underline s^{\eta} : x\mapsto \left\{ 
\begin{array}{lll}
u_1^\star & \textrm{ if } & x< -\eta,\\
1 & \textrm{ if } & x\in (-\eta ,0),\\
u_2^\star & \textrm{ if }& x >0,
\end{array}\right.
$$
and 
$$
\overline s^{\eta} : x\mapsto \left\{ 
\begin{array}{lll}
1 & \textrm{ if } & x< 0,\\
 u_2^\star & \textrm{ if } & x\in (0,\eta),\\
 1 & \textrm{ if }& x >\eta.
\end{array}\right.
$$
\end{proposition}
\vskip 10pt
\begin{proof}
We set 
\be\label{eq:us}
\underline{s}^{\eps,\eta}(x) = \left\{\begin{array}{lll}
y^\eta\left(\frac{x+\eta}{\eps} -\eta \right) & \textrm{ if } & x< 0,\\
u_2^\star & \textrm{ if }& x >0,
\end{array}\right.
\ee
and
 \be\label{eq:os}
\overline{s}^{\eps,\eta}(x) = \left\{\begin{array}{lll}
1 & \textrm{ if }Ê& x<0, \\
z^\eta\left(\frac{x-\eta}{\eps} + \eta \right) & \textrm{ if } & x> 0,
\end{array}\right.
\ee
where the functions $y^\eta$ and $z^\eta$ have been defined in Lemmas~\ref{EDO_lem} and \ref{lem:EDO_2}. Since the functions $\phii(\underline{s}^{\eps,\eta})$ and $\phii(\overline{s}^{\eps,\eta})$ are monotone in $\O_i$, there derivatives $\frac{\rm d}{{\rm d}x} \phii(\underline{s}^{\eps,\eta})$ and  $\frac{\rm d}{{\rm d}x} \phii(\overline{s}^{\eps,\eta})$ belong to $L^1(\R)$, and also to $L^\infty(\R)$ because $\underline{s}^{\eps,\eta}$ and $\overline{s}^{\eps,\eta}$ are solutions to~\eqref{eq:steady}. Thus they belong to $L^2(\R)$. Hence, for fixed $\eps$, 
$\underline{s}^{\eps,\eta}$ and $\overline{s}^{\eps,\eta}$ are bounded flux solutions to the problem~\eqref{Pe_NC}. The convergence 
as $\eps\to0$ towards the functions $\overline s^\eta$ and $\underline s^\eta$ is a direct consequence of Lemmas~\ref{EDO_lem} and \ref{lem:EDO_2}.
\end{proof}

\subsection{a ${ L^2((0,T);H^1(\O_i))}$ estimate}

Our goal is now to derive an estimate which ensures that the effects of capillarity 
vanish almost everywhere in $\O_i\times\R_+$ as $\eps$ tends to $0$.
\vskip 10pt
\begin{proposition}\label{L2H1_NC}Let $u_0\in L^\infty(\R)$ with $0 \le u_0 \le 1$ a.e. satisfying the assumptions of Theorem~\ref{comp_bounded_NC} and let $\ue$ be the corresponding bounded flux solution. Then for all $\eps\in(0,1)$, for all $T>0$,
there exists $C$ depending only on $u_0,g_i,\phii,T$ such that 
\be\label{L2H1_NC_2}
\sqrt{\eps} \| \partial_x \phii(\ue) \|_{L^2(\OiT)} \le C.\ee
This particularly ensures that 
\be\label{L2H1_NC_3}
\eps \| \partial_x\phii(\ue)\|_{L^2(\OiT)} \to 0 \quad  \textrm{ as } \eps\to0.
\ee
\end{proposition}
\vskip 10pt
The idea of the proof of Proposition~\ref{L2H1_NC} is formally to choose 
$(\ue-\hat u)\psi$ as test function in \refe{bounded_for_NC} for a function $x\mapsto\psi(x)$ 
compactly supported in $\O_i$. Using the fact that the flux $F_i^\eps$ is uniformly 
bounded in $L^\infty(\OiT)$, we can let $\psi$ tend towards $\chi_{\O_i}$, with 
$\chi_{\O_i}(x)=1$ if $x\in \O_i$ and $0$ otherwise, and the estimate \refe{L2H1_NC_2} 
follows. To obtain \refe{L2H1_NC_3}, it suffices to multiply \refe{L2H1_NC_2} by 
$\sqrt\eps$. We refer to \cite[Proposition~2.3]{NPCX} for a more details on the proof of 
Proposition~\ref{L2H1_NC}.

\subsection{approximation of the initial data}
In order to ensure that the limit $u$ of the approximate solutions $\ue$ as $\eps\to 0$ admits the expected strong traces on the interface $\{x=0\}$, 
we will perturb the initial data $u_0$. 
\vskip 10pt

\begin{lemma}\label{lem:init_eta_eps}
Let $u_0\in L^\infty(\R)$ satisfying~\eqref{hyp:large_init}, then there exists $\left(u_0^{\eps,\eta}\right)_{\eps,\eta}$ such that 
\begin{enumerate}
\item[(a).] $\underline s^{\eps,\eta} (x) \le u_0^{\eps,\eta}(x) \le \overline s^{\eps,\eta}(x)$ a.e. in $\R$, where the functions $\underline s^{\eps,\eta}$ and $\overline s^{\eps,\eta}$ are defined in~\eqref{eq:us}-\eqref{eq:os},
\vskip 5pt
\item[(b).] $\eps\left\|\partial_x \varphi_i(u_0^{\eps,\eta})\right\|_{L^\infty(\O_i)} \le C$ where $C$ depends neither on $\eps$ nor on $\eta$, 
\vskip 5pt
\item[(c).] $u_0^{\eps,\eta} \to u_0$ in  $L^1_{loc}(\R)$ as $(\eps,\eta)\to (0,0)$.
\end{enumerate}
\end{lemma}
\vskip 10pt
\begin{proof}
Let $\left(\rho_n\right)_{n\in \N^\star}$ be a sequence of mollifiers, then $\rho_n\ast u_0$ is a smooth function tending $u_0$ as $n\to\infty$. 
Then, for $\eps>0$, we choose $n\in \N^\star$ such that 
\be\label{eq:init_45}
\max\left\{ n, \left\| \partial_x \phii(u_0\ast \rho_n) \right\|_{L^\infty(\O_i)} \right\} \ge \frac1\eps,
\ee
and we define
\be\label{eq:inti_eps_eta}
u_0^{\eps,\eta}(x) = \max\left\{ \underline s^{\eps,\eta}(x), \min\left\{ \overline{s}^{\eps,\eta}(x), u_0\ast\rho_n(x)\right\}\right\}.
\ee
The point (a) is a direct consequence of~\eqref{eq:inti_eps_eta}. Letting $(\eps,\eta)\to (0,0)$ in~\eqref{eq:inti_eps_eta} yields 
$$
\lim_{(\eps,\eta)\to(0,0)} u_0^{\eps,\eta}(x) = \max\left\{ u_i^\star, \min\left\{1, u_0(x)\right\}\right\}.
$$
Since $u_0$ is supposed to satisfy~\eqref{hyp:large_init}, this provides 
$$
\lim_{(\eps,\eta)\to(0,0)} u_0^{\eps,\eta}(x) = u_0(x) \quad \textrm{ a.e. in } \R.
$$
The point (c) follows. In order to establish (b), it suffices to note that there exist an open subset $\omega$ of $\R$ such that 
$u_0^{\eps,\eta}(x)$ is equal to $u_0\ast \rho_n(x)$ for $x\in\omega$, and such that 
$u_0^{\eps,\eta}(x)$ is either equal to $\underline s^{\eps,\eta}(x)$ or to $\overline s^{\eps,\eta}(x)$ on $\omega^c = \R\setminus\omega$.
It follows from~\eqref{eq:init_45} that 
$$\eps \left\|\partial_x \phii(u_0^{\eps,\eta}) \right\|_{L^\infty(\O_i\cap\omega)} \le 1.$$
One has 
$$
f_i(u_0^{\eps,\eta})(x) - \eps \partial_x \phii(u_0^{\eps,\eta})(x) = q\quad \textrm{ a.e. in } \O_i\cap\omega^c,
$$
thus
$$\eps \left\|\partial_x \phii(u_0^{\eps,\eta}) \right\|_{L^\infty(\O_i\cap\omega^c)} \le \| q  - f_i \|_{L^\infty(u_i^\star,1)} .$$
This concludes the proof of Lemma~\ref{lem:init_eta_eps}.
\end{proof}
\vskip 10pt
\begin{definition}\label{Def:prepared}
A function $u_0$ is said to be a prepared initial data if it satisfies $(1-u_0)\in L^1(\R)$, $\partial_x \phii(u_0) \in L^\infty(\O_i)$ and 
\be\label{eq:prepared}
\underline s^{\eps,\eta} (x) \le u_0(x) \le \overline s^{\eps,\eta}(x)\textrm{ a.e. in } \R
\ee
for some $\eps>0,\eta>0$. 
\end{definition}
\vskip 10pt

Since the function $(\eps,\eta) \mapsto \underline s^{\eps,\eta}$ is decreasing with respect to both arguments and since the 
function $(\eps,\eta) \mapsto \overline s^{\eps,\eta}$ is increasing with respect to both arguments, if $u_0$ 
satisfies~\eqref{eq:prepared} for $\eps=\eps_0$ and $\eta=\eta_0$, then $u_0$ satisfies~\eqref{eq:prepared} for all 
$(\eps,\eta)$ such that $\eps\le \eps_0$ and $\eta\le\eta_0$. So the following Proposition is a direct consequence from~\eqref{eq:comp_bounded}. 
\vskip 10pt
\begin{proposition}\label{prop:comp_prepared}
Let $u_0$ be a prepared initial data satisfying~\eqref{eq:prepared} for $\eps=\eps_0$ and $\eta = \eta_0$, then for all $\eps\le \eps_0$, 
for all $\eta\le\eta_0$, the solution $\ue$ to \eqref{Pe_NC} satisfies 
$$
\underline{s}^{\eps,\eta}(x) \le \ue(x,t) \le \overline{s}^{\eps,\eta}(x) \quad\textrm{ for a.e. } (x,t) \in \R\times\R_+.
$$
\end{proposition}

%
%
\section{Convergence}\label{conv_NC_sec}
\subsection{a compactness result}
Since ${(\ue)}_\eps$ is uniformly bounded between 0 and 1, there exists $u\in L^\infty(\R\times(0,T))$ 
such that $\ue \to u$ is the $L^\infty$ weak-star sense. This is of course insufficient to pass 
in the limit in the nonlinear terms. Either greater estimates are needed, like for example a 
$BV$-estimate introduced in the work of Vol$'$pert \cite{Vol67} and in \cite{NPCX}, or  
we have to use a weaker compactness result. This idea motivates the introduction 
of Young measures as in the papers of DiPerna \cite{DP85} and Szepessy \cite{Sze91}, 
or equivalently the notion of nonlinear weak star convergence, introduced in \cite{EGGH98} 
and \cite{EGH00}, which leads to the notion of process solution given in 
Definition~\ref{process_def}.
\vskip 10pt
\begin{theorem}[Nonlinear weak star convergence]\label{NLW*}
Let $\Qq$ be a Borelian subset of $\R^k$, and $(u_n)$ be a bounded sequence in $L^\infty(\Qq)$. 
Then there exists $u\in L^\infty(\Qq\times(0,1))$, such that up to a subsequence, $u_n$ tends to $u$ 
"in the non linear weak star sense" as $n\to \infty$, i.e.: $\forall g\in \Cc(\R,\R), $
$$
g(u_n)\to \int_0^1 g(u(\cdot,\a))d\a \textrm{ for the weak star topology of }
L^\infty(\Qq) \textrm{ as }n\to\infty.
$$
\end{theorem}
\vskip 10pt
We refer to \cite{DP85} and \cite{EGH00} for the proof of Theorem~\ref{NLW*}.

\subsection{convergence towards a process solution}
Because of the lack of compactness, we have to introduce the notion of process solution, inspired from the notion of {\em measure valued solution} introduced by DiPerna \cite{DP85}.
\vskip 10pt
\begin{definition}[process solution to \refe{Plim}]\label{process_def}
A function $u\in L^\infty(\R\times\R_+\times(0,1))$ is said to be a process solution to 
\refe{Plim} if $0\le u \le 1$ and for $i=1,2$, $\forall \psi\in \Dd^+(\overline\O_i\times\R_+)$, $\forall\k\in[0,1]$,
\beqn
\lefteqn{\ds\int_{\R_+}\int_{\O_i}\int_0^1 (u(x,t,\a)-\k)^\pm \partial_t\psi(x,t)d\a dxdt + 
\int_{\O_i} (u_0(x)-\k)^\pm \psi(x,0)dx }\nn\\
&&\ds + \int_{\R_+}\int_{\O_i}\int_0^1 \Phi_{i\pm}(u(x,t,\a),\k) \partial_x\psi(x,t) d\a dxdt  
+ M_{f_i}\int_{\R_+} (\overline u_i-\k)^\pm \psi(0,t)dt \ge 0,\nn
\eeqn
where $M_{f_i}$ is any Lipschitz constant of $f_i$, $\overline u_1=1$ and  $\overline u_2=u_2^\star$.
\end{definition}
\vskip 10pt
\begin{lemma}\label{lem:strong_traces}
Let $u_0$ be a $\eta$-prepared initial data in the sense of Definition~\ref{Def:prepared} for some $\eta>0$, and let ${(\ue)}_\eps$ be the corresponding
family of approximate solutions. Then 
\be\label{eq:strong_trace_1}
\ue(x,t) \to 1\quad \textrm{ for a.e. }(x,t)\in (-\eta,0)\times\R_+,
\ee
\be\label{eq:strong_trace_2}
\ue(x,t) \to u_2^\star \quad\textrm{ for a.e. }(x,t)\in (0,\eta)\times\R_+.
\ee
\end{lemma}
\vskip 10pt
\begin{proof}
Firstly, since $u_0$ is a $\eta$-prepared initial data, there exists $\eps_0>0$ such that 
$$
\underline{s}^{\eps_0,\eta} \le u_0 \le \overline{s}^{\eps_0,\eta}.
$$
Then it follows from Proposition~\ref{prop:comp_prepared} that for all $\eps\in(0,\eps_0)$, for a.e. $(x,t)\in \R\times\R_+$
\be\label{eq:comp_47}
\underline{s}^{\eps,\eta}(x) \le \ue(x,t) \le \overline{s}^{\eps,\eta}(x).
\ee
This particularly shows that for all $\eps\in(0,\eps_0)$, for a.e. $(x,t)\in (-\eta,0)\times\R_+$,
$$
\ue(x,t) = 1,
$$
thus \eqref{eq:strong_trace_1} holds. The assertion~\eqref{eq:strong_trace_2} can be obtained by using Proposition~\ref{prop:sub-super} and the dominated convergence theorem.
\end{proof}
\vskip 10pt
\begin{proposition}[convergence towards a process solution]\label{conv_proc}
Let $u_0$ be a prepared initial data in the sense of Definition~\ref{Def:prepared}, and let ${(\ue)}_\eps$ be the corresponding
family of approximate solutions. Then, up to an extraction, $u^\eps$ converges in the nonlinear 
weak-star sense towards a process solution $u$ to the problem \refe{Plim}.
\end{proposition}

\vskip 10pt
\begin{proof}
Since $\ue$ is a weak solution of \refe{Pe_NC}, which is a non-fully degenerate parabolic problem, i.e. 
$\phii^{-1}$ is continuous, it follows from the work of Carrillo \cite{Car99} that $\ue$ is an entropy 
weak solution, i.e.: $\forall \psi\in \Dd^+(\O_i\times\R_+)$, $\forall\k\in[0,1]$,
\beqn
\lefteqn{  \int_{\R_+}\int_{\O_i} (\ue(x,t)-\k)^\pm \partial_t \psi(x,t) dxdt
+ \int_{\O_i} (u_0(x)-\k)^\pm \psi(x,0) dx}\nn\\
&&\ds\hspace{10pt} + \int_{\R_+} \int_{\O_i} \left[\Phi_{i\pm}(\ue(x,t),\k)- \eps\partial_x(\phii(\ue)(x,t)-\phii(\k))^\pm \right] 
\partial_x \psi(x,t)dxdt\ge 0.\nn\label{entro_for_NC}
\eeqn
This family of inequalities is only available for non-negative functions $\psi$ compactly 
supported in $\O_i$, and so vanishing on the interface $\{x=0\}$. To overpass this difficulty, we 
use cut-off functions $\chi_{i,\delta}$. \label{cut-off}

Let $\delta>0$, we denote by $\chi_{i,\delta}$ a smooth non-negative function, 
with $\chi_{i,\delta}(x)=0$ if $x\notin \O_i$, and 
$\chi_{i,\delta}(x)=1$ if $x\in\O_i$, $|x|\ge \delta$.
Let $\psi\in \Dd^+(\overline\O\times\R_+)$, then $\psi\chi_{i,\delta}\in \Dd^+(\O_i\times\R_+)$ can be used 
as test function in \refe{entro_for_NC}. This yields
\beqn
\lefteqn{\ds  \int_{\R_+}\int_{\O_i} (\ue-\k)^\pm \partial_t \psi \chi_{i,\delta}dxdt
+ \int_{\O_i} (u_0-\k)^\pm \psi(\cdot,0)\chi_{i,\delta} dx} \nn\\
&&\ds \hspace{20pt}+ \int_{\R_+} \int_{\O_i} \left[\Phi_{i\pm}(\ue,\k)- \eps\partial_x(\phii(\ue)-\phii(\k))^\pm \right] 
\partial_x \psi \chi_{i,\delta}dxdt \nn \\
&&\ds \hspace{20pt} +\!\!\int_{\R_+} \int_{\O_i}\left[\Phi_{i\pm}(\ue,\k)- \eps\partial_x(\phii(\ue)-\phii(\k))^\pm \right] 
\psi\partial_x  \chi_{i,\delta}dxdt\ge 0.\label{entro_for_NC2}
\eeqn

We can now let $\eps$ tend to $0$. Thanks to Theorem~\ref{NLW*}, there exists $u\in L^\infty(\R\times\R_+\times(0,1))$ such that
\begin{eqnarray}
\lefteqn{ \lim_{\eps\to0}  \int_{\R_+}\int_{\O_i} (\ue(x,t)-\k)^\pm \partial_t \psi(x,t) \chi_{i,\delta}(x)dxdt = } \nn  \\
&&\hspace{50pt} \int_{\R_+}\int_{\O_i} \int_0^1(u(x,t,\a)-\k)^\pm \partial_t \psi(x,t) \chi_{i,\delta}(x)d\a dxdt, \label{eq:NLW*_1}
\end{eqnarray}
\begin{eqnarray}
\lefteqn{ \lim_{\eps\to0}  \int_{\R_+}\int_{\O_i} \Phi_{i\pm}(\ue(x,t),\k) \partial_x \psi(x,t) \chi_{i,\delta}(x)dxdt = } \nn  \\
&&\hspace{50pt} \int_{\R_+}\int_{\O_i} \int_0^1 \Phi_{i\pm}(u(x,t,\a),\k) \partial_x \psi(x,t) \chi_{i,\delta}(x)d\a dxdt . \label{eq:NLW*_2}
\end{eqnarray}
Thanks to Proposition~\ref{L2H1_NC}, one has
$$
\eps\partial_x(\phii(\ue)-\phii(\k))^\pm \textrm{ tends to } 0 \quad\textrm{ a.e. in } \OiT \textrm{ as }\eps\to 0,
$$
then
\be\label{lim_phi(u)}
\lim_{\eps\to0}  \int_{\R_+}\int_{\O_i} \eps\partial_x(\phii(\ue)(x,t) - \phii(\k))^\pm \partial_x \left( \psi(x,t) \chi_{i,\delta}(x)\right) dxdt = 0.
\ee 
Since $u_0$ is supposed to be a $\eta$-prepared initial data for some $\eta>0$, we can claim thanks to 
Lemma~\ref{lem:strong_traces} that $\ue(x,t)$ converges almost everywhere on $(-\eta,\eta)\times\R_+$ towards $\overline u_i$ if 
$x\in\O_i$. Since for $\delta<\eta$ small enough, the support of $\partial_x \chi_{1,\delta}$ is 
included in the set where $\ue$ converges strongly, one has
\begin{eqnarray}
\lefteqn{ 
\lim_{\eps\to0}  \int_{\R_+}\int_{\O_i} \Phi_{i\pm}(\ue(x,t),\k)  \psi(x,t) \partial_x\chi_{i,\delta}(x)dxdt = } \nn \\
&&\hspace{100pt} 
\int_{\R_+}\int_{\O_i} \Phi_{i\pm}(\overline u_i,\k) \psi(x,t) \partial_x \chi_{i,\delta}(x)dxdt . \label{eq:NLW*_3}
\end{eqnarray}
We let now $\delta$ tend to $0$. Since $\chi_{i,\delta}(x)$ tends to $1$ a.e. in $\O_i$, \eqref{eq:NLW*_1} and \eqref{eq:NLW*_2} respectively provide
\begin{eqnarray}
\lefteqn{ \lim_{\delta\to0} \lim_{\eps\to0}  \int_{\R_+}\int_{\O_i} (\ue(x,t)-\k)^\pm \partial_t \psi(x,t) \chi_{i,\delta}(x)dxdt = } \nn  \\
&&\hspace{100pt} \int_{\R_+}\int_{\O_i} \int_0^1(u(x,t,\a)-\k)^\pm \partial_t \psi(x,t) d\a dxdt, \label{eq:NLW*_4}
\end{eqnarray}
\begin{eqnarray}
\lefteqn{ \lim_{\delta\to 0}\lim_{\eps\to0}  \int_{\R_+}\int_{\O_i} \Phi_{i\pm}(\ue(x,t),\k) \partial_x \psi(x,t) \chi_{i,\delta}(x)dxdt = } \nn  \\
&&\hspace{100pt} \int_{\R_+}\int_{\O_i} \int_0^1 \Phi_{i\pm}(u(x,t,\a),\k) \partial_x \psi(x,t) d\a dxdt . \label{eq:NLW*_5}
\end{eqnarray}
One has also 
\be\label{eq:NLW*_6}
\lim_{\delta\to 0} \int_{\O_i} (u_0(x)-\k)^\pm \psi(x,0)\chi_{i,\delta}(x) dx = \int_{\O_i} (u_0(x)-\k)^\pm \psi(x,0)dx.
\ee
One has $$\left| \Phi_{i\pm}(\overline u_i,\k) \right| \le M_{f_i} \left( \overline u_i - \k \right)^\pm$$
then
\begin{eqnarray*}
\lefteqn{\left|\int_{\R_+}\int_{\O_i} \Phi_{i\pm}(\overline u_i,\k) \psi(x,t) \partial_x \chi_{i,\delta}(x)dxdt\right|}\\
&&\hspace{100pt}\le M_{f_i} \left( \overline u_i - \k \right)^\pm \int_{\R_+}\int_{\O_i}  \psi(x,t)\left|  \partial_x \chi_{i,\delta}(x) \right| dxdt.
\end{eqnarray*}
Since $\left| \partial_x \chi_{i,\delta} \right|$ tends to $\delta_{x=0}$ in the $\Mm(\R)$-weak star sense where $$\Big<\delta_{x=0},\zeta\Big>_{\Mm(\R),\Cc_0(\R)} = \zeta(0),$$
we obtain that 
\begin{eqnarray}
\lefteqn{ 
\liminf_{\delta\to0}\lim_{\eps\to0}  \int_{\R_+}\int_{\O_i} \Phi_{i\pm}(\ue(x,t),\k)  \psi(x,t) \partial_x\chi_{i,\delta}(x)dxdt \ge 
} \nn \\
&&\hspace{180pt} 
M_{f_i} \left( \overline u_i - \k \right)^\pm\int_{\R_+}\psi(0,t)dt . \label{eq:NLW*_7}
\end{eqnarray}
Using~\eqref{lim_phi(u)},\eqref{eq:NLW*_4},\eqref{eq:NLW*_5},\eqref{eq:NLW*_6},\eqref{eq:NLW*_7} in \eqref{entro_for_NC2} shows that $u$ is a process solution in the sense of Definition~\ref{process_def}.
\end{proof}

\subsection{uniqueness of the (process) solution}

It is clear that the notion of process solution is weaker than the one of solution
given in Definition~\ref{Sol_NC_def}.
We state here a theorem which claims the equivalence of the two notions, 
i.e. any  process solution is a solution in the sense of Definition~\ref{Sol_NC_def}. 
Furthermore, such a solution is unique, and a $L^1$-contraction principle can 
be proven.
\vskip 10pt
\begin{theorem}[uniqueness of the (process) solutions]\label{uni_NC}
There exists a unique process solution $u$ to the problem~\refe{Plim}, and 
furthermore this solution does not depend on $\a$, i.e. $u$ is a solution 
to the problem \refe{Plim} in the sense of definition~\ref{Sol_NC_def}.
Furthermore, if $u_0$, $v_0$ are two initial data in $L^\infty(\R)$ satisfying~\eqref{hyp:large_init} and let
$u$ and $v$ be two solutions associated to those initial data, then
for all $t\in[0,T)$,
\be\label{contract_NC}
\int_{-R}^{R} (u(x,t)-v(x,t))^\pm dx \le \int_{-R - M_f t}^{R+M_f t} (u_0(x)-v_0(x))^\pm dx.
\ee
\end{theorem}
\vskip 10pt
This theorem is a consequence of \cite[Theorem 2]{Vov02}. Let $u(x,t,\a)$ and $v(x,t,\b)$ be two process solutions corresponding to initial data $u_0$ and $v_0$. Classical Kato inequalities can be derived in each $\O_i \times\R_+$ by using the doubling variable technique: $\forall \psi\in \Dd^+(\O_i\times\R_+)$, 
\begin{eqnarray*}
\lefteqn{\int_{\R_+}\int_{\O_i}\int_0^1\int_0^1 (u(x,t,\a) - v(x,t,\b))^\pm\partial_t \psi(x,t) d\a d\b dxdt }\\
&&\hspace{50pt}+ \int_{\O_i} (u_0(x) - v_0(x))^\pm \psi(x,0) dx \\
&&\hspace{50pt}+ \int_{\R_+}\int_{\O_i}\int_0^1\int_0^1\Phi_{i\pm}(u(x,t,\a),v(x,t,\b)) \partial_x \psi(x,t) d\a d\b dxdt \ge 0.
\end{eqnarray*}
The treatment of the boundary condition at the interface is an adaptation to the case of process solution 
to the work of Otto summarized in~\cite{Otto96_CRAS} and detailed in~\cite{MNRR96} leading to (see 
\cite[Lemma 2]{Vov02}): $\forall \psi\in \Dd^+(\overline\O_i\times\R_+)$, 
\begin{eqnarray}
\lefteqn{\int_{\R_+}\int_{\O_i}\int_0^1\int_0^1 (u(x,t,\a) - v(x,t,\b))^\pm\partial_t \psi(x,t) d\a d\b dxdt } \nn \\
&&\hspace{30pt}+ \int_{\O_i} (u_0(x) - v_0(x))^\pm \psi(x,0) dx \nn\\
&&\hspace{30pt}+ \int_{\R_+}\int_{\O_i}\int_0^1\int_0^1\Phi_{i\pm}(u(x,t,\a),v(x,t,\b)) \partial_x \psi(x,t) d\a d\b dxdt \ge 0. \label{eq:Kato_2}
\end{eqnarray}
Choosing 
$$\psi_\eps(x,s) = \left\{\begin{array}{lll}1&\textrm{ if } &|x|\le R+ M_f s ,\\
\ds \frac{R+ M_f s + \eps - |x|}{\eps}&\textrm{ if } & R+ M_f t \le |x| \le R+ M_f s + \eps \\
0  &\textrm{ if } &|x|\ge R+ M_f s + \eps
\end{array}\right.$$
if $s \le t$ and $\psi_{\eps}(x,s)=0$ if $s > t$ 
as test function in~\eqref{eq:Kato_2} and letting $\eps$ tend to $0$ provide the expected $L^1$-contraction principle~\eqref{contract_NC}.
\vskip 10pt
Finally, if $u$ and $\t u$ are two process solutions associated 
to the same initial data $u_0$, we obtain a $L^1$-contraction principle of the 
following form: for a.e. $t\in\R_+$,
$$
\int_\R \int_0^1 \int_0^1 (u(x,t,\a)-\t u(x,t,\beta))^\pm d\a d\beta dx \le 0,
$$
thus $u(x,t,\a) = \t u(x,t,\b)$ a.e. in $\R\times\R_+\times(0,1)\times(0,1)$. Hence $u$ does not depend on the process variable $\a$.
\vskip 10pt
\begin{theorem}\label{conv_en_NC}
Let $u_0$ be a prepared initial data in the sense of Definition~\ref{Def:prepared}, and let $\ue$ be the corresponding solution to 
the approximate problem \refe{Pe_NC}. Then $\ue$ converges to the unique 
solution $u$ to \eqref{Plim} associated to initial data $u_0$ in the $L^p((0,T);L^q(\R))$-sense, for all 
$p,q\in [1,\infty)$.
\end{theorem}
\vskip 10pt
\begin{proof}
We have seen in Proposition~\ref{conv_proc} that $\ue$ converges up to an extraction 
towards a process solution. The family ${(\ue)}_\eps$  admits so a unique 
adherence value, which is a solution thanks to Theorem~\ref{uni_NC}, 
thus the whole family converges towards this unique limit $u$. 
\vskip 10pt
Let $K$ denotes a compact subset of $\R\times[0,T]$, then one has
$$
\iint_K (\ue - u)^2 dxdt = \iint_K \left(\ue\right)^2 dx -2 \iint_K \ue u dx + \iint_K u^2 dx.
$$
Since $\ue$ converges in the nonlinear weak star sense towards $u$,  
$$
\lim_{\eps\to 0} \iint_K \left(\ue\right)^2 dx =  \iint_K u^2 dx.
$$
Moreover, $\ue$ converges in the $L^\infty$ weak star topology towards $u$, 
then 
$$
\lim_{\eps\to 0} \iint_K \ue u dx =  \iint_K u^2 dx.
$$
Thus we obtain
$$
\lim_{\eps\to 0}\iint_K (\ue - u)^2 dxdt =0.
$$\label{ae_conv}
One concludes using the fact the $|\ue - u| \le 1$ for all $\eps>0$.
\end{proof}

\subsection{initial data in $L^\infty(\R)$}
In this section, we extend the result of Theorem~\ref{conv_en_NC} to any initial data in $L^\infty(\R)$ satisfying~\eqref{hyp:large_init} thanks to density argument.
 
\vskip 10pt
\begin{theorem}
Let $u_0\in L^\infty(\R)$ satisfying \eqref{hyp:large_init}, and let $\left(u_{0,n}\right)_{n\in \N^\star}$ be a sequence of prepared initial data tending to $u_0$ in $L^1_{loc}(\R)$. Then the sequence $\left(u_n\right)_n$ of solutions to~\eqref{Plim} corresponding to the sequence $\left(u_{0,n}\right)$ of initial data converges in $\Cc(\R_+;L^1_{loc}(\R))$ towards the unique solution to~\eqref{Plim} corresponding to solution the initial data $u_0$.
 \end{theorem}
 \vskip 10pt
 \begin{proof}
 First, note that for all $u_0\in L^\infty(\R)$ satisfying \eqref{hyp:large_init}, there exists a sequence $\left(u_{0,n}\right)_{n\in \N^\star}$ of prepared initial data tending to $u_0$ in $L^1_{loc}(\R)$ thanks to Lemma~\ref{lem:init_eta_eps}. 
\vskip 10pt
 Thanks to~\eqref{contract_NC}, one has for $n,m\in \N^\star$, for all $t\in\R_+$
 $$
 \int_{-R}^{R} (u_n(x,t)-u_m(x,t))^\pm dx \le \int_{-R - M_f t}^{R+M_f t} (u_{0,n}(x)-u_{0,m}(x))^\pm dx,
 $$
 then $\left(u_n\right)_n$ is a Cauchy sequence in $\Cc(\R_+;L^1_{loc}(\R))$. In particular, there exists $u$ such that 
 $$
 \lim_{n\to\infty} u_n = u \quad \textrm{ in } \Cc(\R_+; L^1_{loc}(\R)).
 $$
 It is then easy to check that $u$ is the unique solution to~\eqref{Plim}.
 \end{proof}

%
%

\section{Entropy solution for small initial data}\label{entro_small}
In this section, we suppose that the initial data $u_0$ belongs to $L^1(\R)$, and that 
\be\label{H_u0_entro}
0 \le u_0 \le u_i^\star\quad \textrm{ a.e. in }\O_i.
\ee
This initial data can be smoothed using following lemma whose proof is almost the same as the proof of Lemma~\ref{lem:init_eta_eps}.
\vskip 10pt
\begin{lemma}\label{lem:prepared_2}
There exists $\left(u_0^\eps\right)_{\eps>0}\subset L^1(\R)$ such that 
\begin{itemize}
\item $\partial_x \phii(u_0^\eps)\in L^\infty(\O_i)$,
\item ${\rm ess}\lim_{x\nearrow 0} u_0^\eps(x) = 1$,
\item $\lim_{\eps\to0} u_0^\eps = u_0$ in $L^1_{loc}(\R)$.
\end{itemize}
\end{lemma}
\vskip 10pt
For all $\eps>0$, there exists a unique bounded flux solution $\ue$ to~\eqref{Pe_NC} corresponding to $u_0^\eps$ thanks to Theorem~\ref{comp_bounded_NC}. The following theorem claims that as $\eps$ tends to $0$, $\ue$ tends to the unique entropy solution in the sense of Definition~\ref{Def:entro}.
\vskip 10pt
\begin{theorem}[convergence towards the entropy solution]
Let $u_0\in L^\infty(\R)$ satisfying~\eqref{H_u0_entro} and let $\left(u_0^\eps\right)_\eps$ be a family of approximate initial data built in Lemma~\ref{lem:prepared_2}. Let $\ue$ be the bounded flux solution to~\eqref{Pe_NC} corresponding to $u_0^\eps$, then $\ue$ converges to $u$ in $L^1_{loc}(\R\times\R_+)$ as $\eps$ tends to $0$ where $u$ is the unique entropy solution in the sense of Definition~\ref{Def:entro}.
\end{theorem}
\vskip 10pt
\begin{proof}
Using the technics introduced in \cite[Proposition 2.8]{NPCX}, we can show that for all $\lambda\in[0,q]$ there exists a steady solution $\k_\lambda^\eps$ to the problem \eqref{Pe_NC}, corresponding to a constant flux 
$$f_i(\k^\eps_\lambda)-\eps \partial \phii(\k^\eps_\lambda) = \lambda,$$
and such that this solution converges uniformly on each compact subset of $\R^\star$ as $\eps$ tends to $0$ towards 
$$\k_\lambda(x)= \min_\k\left\{f(\k,x) = \lambda\right\}.$$ 
Following the idea of Audusse and Perthame~\cite{AP05}, we will now compare the limit $u$ of $\ue$ as $\eps$ to $0$ with the steady state $\k_\lambda$. Let $\l\in [0,q]$. Since $\ue$ and $\k^\eps_\l$ are both bounded flux solutions, it follows from Proposition~\ref{comp_prop} that for all $\psi\in \Dd^+(\R\times{\R_+})$,
\beqn
\lefteqn{\ds \int_{\R_+} \int_\R \left( \ue -  \k^\eps_\l\right)^\pm \partial_t \psi dxdt + \int_\R \left( u_0^\eps -  \k^\eps_\l\right)^\pm \psi(\cdot,0) dx}\nn\\
&&\ds\hspace{20pt} + \int_{\R_+} \sum_i \int_{\O_i} \left(  \Phi_{i\pm}(\ue, \k_\l^\eps) - \eps \partial_x \left( \phii(\ue) - \phii(\k_\l^\eps)\right)^+ \right)  \partial_x\psi dxdt \ge 0. \label{comp_kle}
\eeqn
Choosing $\l=q$ and $\psi(x,t)=(T-t)^+ \xi(x)$ for some arbitrary $T>0$ and some $\xi\in\Dd^+(\R)$ yields
\be\int_0^T \int_\O (\ue-\k^\eps_q)^+  \xi dxdt \le \int_0^T (T-t) \sum_{i=1,2} \int_{\O_i} 
\eps \partial_x \left(\phii(\ue) - \phii(\k^\eps_q)\right)^+ \partial_x \xi dxdt. \label{comp_ke}
\ee
Since $\ue$ is bounded between $0$ and $1$, it converges in the nonlinear weak star sense, thanks to Theorem~\ref{NLW*} towards a function $u\in L^\infty(\R\times{\R_+}\times(0,1))$, with 
$0 \le u \le 1$ a.e.. Then \eqref{comp_ke} provides 
\be
u \le \k_q = u_i^\star \quad \textrm{ a.e. in }\O_i\times{\R_+}\times(0,1).  \label{max_entro}
\ee
Let $\l\in [0,q]$, then taking the limit for $\eps\to 0$ in \eqref{comp_kle} yields
\beqn
\lefteqn{\ds \int_{\R_+} \int_\R\int_0^1 \left| u -  \k_\l\right| \partial_t \psi d\a dxdt + \int_\R \left| u_0 -  \k_\l\right| \psi(\cdot,0) dx}\nn\\
&&\ds\hspace{100pt} + \int_{\R_+} \sum_i \int_{\O_i} \int_0^1 \Phi_{i}(u, \k_\l)  \partial_x\psi d\a dxdt \ge 0. \label{comp_kl}
\eeqn
Suppose that $u_2^\star\ge u_1^\star$. Let $\k\in [0,u_2^\star]$, we denote by $\tilde\k = f_1^{-1}(f_2(\k))\cap [0,u_1^\star]$. Then choosing $\l = f_2(\k)$ in \eqref{comp_kl}, and letting $\eps$ tend to $0$ gives: $\forall\k\in [0,u_2^\star]$, $\forall \psi\in\Dd^+(\R\times{\R_+})$, 
\begin{eqnarray}
\lefteqn{\int_0^T \int_{\O_1} \int_0^1 | u - \t \k| \partial_t \psi d\a dxdt+ \int_{\O_1} |u_0- \t \k | \psi(\cdot,0) dx}\nn\\
 &&\hspace{50pt}+ \ds \int_{\R_+} \int_{\O_2} \int_0^1 | u - \k| \partial_t \psi d\a dxdt+ \int_{\O_2} |u_0-\k | \psi(\cdot,0) dx \nn\\
&&\hspace{50pt}+ \int_{\R_+} \int_0^1 \left( \int_{\O_1} \Phi_1(u, \t \k)  \partial_x \psi dx + \int_{\O_2} \Phi_2(u,\k)  \partial_x \psi dx\right)d\a dt \ge 0. \label{comp2_k}
\end{eqnarray}
It follows from the work of Jose Carrillo \cite{Car99} that the following entropy inequalities hold for test functions compactly supported in $\O_1$: $\forall\k\in [0,1]$, $\forall \psi\in\Dd^+(\O_1\times{\R_+})$,
\begin{eqnarray}
 \lefteqn{ \ds  \int_{\R_+} \int_{\O_1} | \ue - \k| \partial_t \psi dxdt+ \int_{\O_1} |u_0^\eps-\k | \psi(\cdot,0) dx        }\nn\\
 &&\hspace{50pt}+ \int_{\R_+} \int_{\O_1} \left( \Phi_1(\ue,\k) -\eps  \partial_x \left| \varphi_1(\ue) - \varphi_1(\k) \right| \right) \partial_x \psi dxdt \ge 0. \label{eq:entro_car}
\end{eqnarray}
Thus letting $\eps$ tend to $0$ in~\eqref{eq:entro_car} provides: $\forall \psi\in\Dd^+(\O_1\times{\R_+})$, $\forall\k\in[0,1]$, 
\begin{eqnarray}
 \lefteqn{\ds  \int_{\R_+} \int_{\O_1} \int_0^1 | u - \k| \partial_t \psi d\a dxdt+ \int_{\O_1} |u_0-\k | \psi(\cdot,0) dx  }\nn\\
 &&\hspace{100pt} + \int_{\R_+} \int_{\O_1} \int_0^1  \Phi_1(u,\k)  \partial_x \psi d\a dxdt \ge 0. \label{comp1_k}
\end{eqnarray}
Let $\delta>0$, and let $\psi\in \Dd^+(\R\times{\R_+})$, we define 
$$\psi_{1,\delta}(x,t) = \psi (x,t) \chi_{1,\delta}(x),\qquad \psi_{2,\delta}=\psi-\psi_{1,\delta},$$ where $\chi_{1,\delta}$ is the cut-off function introduced in section~\ref{cut-off}. Then using $\psi_{1,\delta}$ as test function in \eqref{comp1_k} and $\psi_{2,\delta}$ in \eqref{comp2_k} leads to:
\beqn
\lefteqn{ \ds \int_{\R_+}\int_\R\int_0^1 | u - \k | \partial_x \psi d\a dxdt + \int_{\R} |u_0-\k | \psi(\cdot,0) dx }\nn \\
&&\ds\hspace{20pt} + \int_{\R_+}\sum_i\int_{\O_i}\int_0^1  \Phi_i(u,\k)  \partial_x \psi d\a dxdt  \nn \\
&&\ds\hspace{20pt}  + \int_{\R_+}\int_{\O_1}\int_0^1 \left( \Phi_1(u,\k) - \Phi_1(u,\t\k) \right) \psi \partial_x \chi_{1,\delta} d\a dxdt\ge   \Rr(\k,\psi,\delta),  \label{comp3_k}
\eeqn
where $\lim_{\delta\to0} \Rr(\k,\psi,\delta) = 0.$ Since $f_1$ is increasing on $[0,u_1^\star]$ and $f_1([u_1^\star,1)) \subset [q,\infty)$, either 
$\k\le u_1^\star$, or $f_1(\k) \ge f_1(u_1^\star).$ This ensures that 
$$\Phi_1(u,\k) = | f_1(u) - f_1(\k) |, \qquad \forall u \in [0,u_1^\star],\quad \forall k\in [0,u_2^\star].$$
This yields
\begin{eqnarray}
\left|\Phi_1(u,\k) - \Phi_1(u,\t\k) \right|&=& \big|| f_1(u) - f_1(\k) | - | f_1(u)- f_1(\t\k) | \big| \nn\\
&\le& | f_1(\k) - f_1(\t\k) | =  | f_1(\k) - f_2(\k) |. \label{triangle_k}
\end{eqnarray}
Taking the inequality \eqref{triangle_k} into account in \eqref{comp3_k}, and letting $\delta\to0$ provides: \\ 
$\forall \k\in [0,\| u \|_\infty]$, $\forall \psi\in \Dd^+(\R\times{\R_+})$, 
\begin{eqnarray*}
\lefteqn{ \int_{\R_+}\int_\R\int_0^1 | u - \k | \partial_x \psi d\a dxdt + \int_{\R} |u_0-\k | \psi(\cdot,0) dx  }  \nn \\
&& \ds\hspace{20pt}  + \int_{\R_+}\sum_i\int_{\O_i}\int_0^1  \Phi_i(u,\k)  \partial_x \psi d\a dxdt 
+ |f_1(\k) - f_2(\k)| \int_0^T \psi (0,\cdot) dt \ge 0.  
\end{eqnarray*}
Using the work of Florence Bachmann~\cite[Theorem 4.3]{Bachmann_These}, we can claim that $u$ is the unique entropy solution to the problem. Particularly, $u$ does not depend on $\a$ (introduced for the nonlinear weak star convergence). As proven in the proof of Theorem~\ref{conv_en_NC}, this implies that $\ue$ converges in $L^1_{loc}(\R\times\R_+)$ towards $u$.
\end{proof}

%
%

\section{Resolution of the Riemann problem}\label{sec:Riemann}
In this section, we complete the resolution of the Riemann problem at the interface $\{x=0\}$, whose result has been given in section~\ref{subsec:oil-trapping}. Consider the initial data 
$$u_0(x) = \left\{Ê\begin{array}{lll}
u_\ell & \textrm{ if } & x<0,\\
u_r & \textrm{ if }Ê& x>0.
\end{array}\right.$$
We aim to determine the traces $(u_1,u_2)$ at the interface  of the solution $u(x,t)$ corresponding to $u_0$. This resolution has already been performed in the following cases.
\begin{itemize}
\item[(a).] $u_1^\star < u_\ell \le 1$ and $u_2^\star \le u_r < 1$: it has been seen that $u_1 = 1$ and $u_2 = u_2^\star$. \vskip 5pt
\item[(b).] $0 \le u_\ell \le u_1^\star$ and $0 \le u_r \le u_2^\star$: Since $u$ is the unique optimal entropy solution studied in \cite{AMV05,Kaa99}, then $u_1 = u_\ell$ and $u_2$ is the unique value in $[0,u_2^\star]$ such that $f_1(u_\ell) = f_2(u_2)$.
\end{itemize}
\vskip 10pt
\noindent 
In the cases
\begin{itemize}
\item[(c).] $u_1^\star < u_\ell \le 1$ and  $u_r = 1$,\vskip 5pt
\item[(d).] $u_\ell = u_1^\star$ and $u_r  = 1$,
\end{itemize}
it is possible to approach the solution $u$ by bounded flux solutions $\ue$ that are constant equal to $1$ in $\O_2\times\R_+$. Then one obtains $u_1 = u_2 = 1$ for the case (c) and $u_1= u_1^\star$  and $u_2 = 1$ for the case (d).
\vskip 10pt
\noindent The last points we have to consider are 
\begin{itemize}
\item[(e).] $u_1^\star < u_\ell \le 1$ and $0 \le u_r < u_2^\star$,
\vskip 5pt 
\item[(f).] $0 \le u_\ell \le u_1^\star$ and $u_2^\star < u_r \le 1$. 
\end{itemize}
To perform the study of the two last cases (e) and (f), we need the following lemmas that can be proved using similar arguments than those used in \cite{FVbarriere}, particularly concerning the treatment of the boundary condition imposed on $\{x=0\}$.
\vskip 10pt
\begin{lemma}\label{lem:e}
Let $u_r\in [0,u_2^\star)$. For all $\eps>0$,  there exists a function $v^\eps$ solution to the problem
\be\label{syst:Pe_2}
\left\{\begin{array}{lll}
\partial_t v^{\eps} + \partial_x \big( f_2(v^{\eps}) - \eps \partial_x \varphi_2(v^{\eps}) \big) = 0 & \textrm{ if } & x>0,\ t>0, \\[5pt]
f_2(v^{\eps}) - \eps \partial_x \varphi_2(v^{\eps}) = f_2(u_2^\star) & \textrm{ if } & x=0,\ t>0, \\[5pt]
v^{\eps} = u_r & \textrm{ if } & x>0 , \ t=0,
\end{array}
\right. 
\ee
satisfying furthermore $u_r\le v^{\eps} \le u_2^\star$ and $\partial_x \varphi_2(v^{\eps}) \in L^\infty(\R_+\times\R_+)$. 
\end{lemma}
\vskip 10pt
\begin{lemma}\label{lem:f}
Let $u_\ell \in [0,u_1^\star]$, $u_r\in (u_2^\star,1]$ and let $u_2$ be the unique value of $[0,u_2^\star]$ such that $f_2(u_2) = f_1(u_\ell)$. 
For all $\eps>0$ there exists a function $w^{\eps}$ solution to the problem 
\be\label{syst:Pe_22}
\left\{\begin{array}{lll}
\partial_t w^{\eps} + \partial_x \big( f_2(w^{\eps}) - \eps \partial_x \varphi_2(w^{\eps}) \big) = 0 & \textrm{ if } & x>0,\ t>0, \\[5pt]
f_2(w^{\eps}) - \eps \partial_x \varphi_2(w^{\eps}) = f_2(u_2) = f_1(u_\ell) & \textrm{ if } & x=0,\ t>0, \\[5pt]
w^{\eps} = u_r & \textrm{ if } & x>0 , \ t=0,
\end{array}
\right. 
\ee
satisfying furthermore $u_2 \le w^{\eps} \le u_r$ and $\partial_x \varphi_2(w^{\eps}) \in L^\infty(\R_+\times\R_+)$.
\end{lemma}
\vskip 10pt
\textbf{The case (e).} Assume that $u_\ell> u_1^\star$ and $u_r< u_2^\star$.
Let $\left(u_0^\eta\right)_\eta$ be a family of initial data such that $\partial_x \phii(u_0^\eta) \in L^\infty(\O_i)$, $u_0^\eta(x) = 1$ for $x \in (-\eta,0)$, $u_0^\eta(x) \in [u_\ell, 1]$ for a.e. $x\in \O_1$, $u_0^\eta(x) = u_r$ a.e. in $\O_2$ and such that
$$
\left\| u_0^\eta - u_\ell \right\|_{L^1(\O_1)} \le 2\eta.
$$
Then thanks to Theorem~\ref{comp_bounded_NC}, there exists a unique bounded flux solution $u^{\eps,\eta}$ to the problem~\eqref{Pe_NC} corresponding to the initial data $u_0^\eta$. 
It is easy to check that the solution defined in $\Omega_2\times\R_+$ by the function $v^\eps$ introduced in Lemma~\ref{lem:e} and coinciding in  $\O_1\times\R_+$ with the unique bounded flux solution corresponding to the initial data 
$$
\tilde u_0^\eta(x) = \left\{Ê\begin{array}{lll}
u_0^\eta(x) & \textrm{ if }Ê& x<0,\\[5pt]
1 & \textrm{ if }Ê& x>0.
\end{array}\right. 
$$
In particular, as $\eps$ tends to $0$, it follows from arguments similar to those developed in the previous sections that $u^{\eps,\eta}$ converges in $L^1_{loc}(\overline\O_i\times\R_+)$ towards the unique entropy solution to the problem
problem 
\be\label{eq:Riemann_e1}
\left\{ \begin{array}{lll}
\partial_t u^\eta + \partial_x f_1(u^\eta) = 0 & \textrm{ if } x<0,\ t>0,\\[5pt]
u^\eta = 1 & \textrm{ if } x=0,\ t>0,\\[5pt]
u^\eta = u_0^\eta & \textrm{ if } x<0,\ t=0.
\end{array}\right.
\ee
\be\label{eq:Riemann_e2}
\left\{ \begin{array}{lll}
\partial_t u + \partial_x f_2(u) = 0 & \textrm{ if } x>0,\ t>0,\\[5pt]
u = u_2^\star & \textrm{ if } x=0,\ t>0,\\[5pt]
u = u_r & \textrm{ if } x>0,\ t=0.
\end{array}\right.
\ee
Note that the trace condition on the interface $\{x=0\}$ in \eqref{eq:Riemann_e2} is fulfilled in a strong sense since $u_r \le u(x,t) \le u_2^\star$ a.e. in $\O_2\times\R_+$ and $f_2$ is increasing on $[u_r, u_2^\star]$. 
\vskip 10pt
The solution to~\eqref{eq:Riemann_e1} depends continuously on the initial data in $L^1_{loc}$. Hence, letting $\eta$ tend to $0$ in~\eqref{eq:Riemann_e1} provides that the limit $u$ of $u^\eta$ is the unique entropy solution to the problem
$$
\left\{ \begin{array}{lll}
\partial_t u+ \partial_x f_1(u) = 0 & \textrm{ if } x<0,\ t>0,\\[5pt]
u = 1 & \textrm{ if } x=0,\ t>0,\\[5pt]
u = u_\ell & \textrm{ if } x<0,\ t=0.
\end{array}\right.
$$
Note that since $u_1^\star < u_\ell \le u  \le 1$ and $\min_{s\in[u,1]} f_1(s) = f_1(1) = q$, the trace prescribed on the interface $\{x= 0\}$ is fulfilled in a strong sense. This particularly yields that in the case (e), the solution to the Riemann problem is given by 
$$
u_1 = 1, \qquad u_2 =u_2^\star.
$$
\vskip 5pt
\textbf{The case (f).} Following the technique used in \cite{NPCX} and in Section~\ref{entro_small}, there exists a unique function 
$u_\ell^\eps$ solution to the problem:
$$
\left\{
\begin{array}{lll}
\ds f_1(u^\eps_\ell)- \eps\frac{\rm d}{{\rm d}x} \varphi_1(u_\ell^\eps) =  f_1(u_\ell) & \textrm{ if } & x<0,\\[5pt]
u_{\ell}^{\eps}(0) = 1 & \textrm{ if } & x=0.
\end{array}\right.
$$
Let $\ue$ be the function defined by 
$$
\ue(x,t) = \left\{ \begin{array}{lll}
u_\ell^\eps(x) & \textrm{Êif }Ê& x<0,\ t\ge0,\\[5pt]
w^\eps(x,t) & \textrm{ if } & x>0,\ t\ge 0,
\end{array}\right.
$$
where $w^\eps$ is the function introduced in Lemma~\ref{lem:f}. 
Then $\ue$ is a bounded flux solution to the problem~\eqref{Pe_NC} in the sense of Definition~\ref{bounded_Def_NC}.
\vskip 10pt
One has 
$$u_\ell^\eps \to u_\ell\quad \textrm{Êin }ÊL^1_{loc}(\O_1)\quad \textrm{Êas } \eps \to 0,$$
and 
$$
w^\eps \to w\quad \textrm{Êin }ÊL^1_{loc}(\O_2\times\R_+)\quad \textrm{Êas } \eps \to 0
$$
where $w$ is the unique solution to 
$$
\left\{
\begin{array}{lll}
\partial_t w + \partial_x f_2(w) = 0 & \textrm{ if } x >0,\ t>0, \\[5pt]
w = u_2 = f_2^{-1}\circ f_1(u_\ell) & \textrm{ if } x=0, \ t>0, \\[5pt]
w= u_r & \textrm{ if } x>0, \ t=0.
\end{array}\right.
$$
Since $w(x,t) \in [u_2,u_r]$ a.e. in $\O_2 \times \R_+$ and since $\min_{s\in [u_2,w]} f_2(s) = f_2(u_2) = f_1(u_\ell)$, the trace $w=u_2$ is satisfied in 
a strong sense on $\{x=0\}$. 
This yields that the solution to the Riemann problem in the case (f) is given by
$$
u_1 = u_\ell, \qquad u_2 = f_2^{-1}\circ f_1(u_\ell).
$$

%
%

\section{Conclusion}
The model presented here shows that for two-phase flows in heterogeneous porous media with negligible dependance of the capillary pressure with respect to the saturation, the good notion of solution is not always the entropy solution presented for example in \cite{AJV04,Bachmann_These}, and  particular care as to be taken with respect to the orientation of the gravity forces. Indeed, some non classical shock can appear at the discontinuities of the capillary pressure field, leading to the phenomenon of oil trapping. We stress the fact that the non classical shocks appearing in our case have a different origin, and a different behavior of those suggested in the recent paper \cite{vDPP07} (see also \cite{LeFloch02}). Indeed, in this latter paper, this lack of entropy was caused by the introduction of the dynamical capillary pressure \cite{HG90,HG93,Pav89}, i.e. the capillary pressure is supposed to depend also on $\partial_t u$. In our problem, the lack of entropy comes only from the discontinuity of the porous medium.
\vskip 5pt
In order to conclude this paper, we just want to stress that this model of piecewise constant capillary pressure curves can not lead to some interesting phenomenon. Indeed, if the capillary pressure functions $\pi_i$ are such that $\pi_1((0,1))\cap \pi_2((0,1))\neq\emptyset$, it appears in \cite[Section 6]{FVbarriere} (see also \cite{BPvD03}) that some oil can overpass the boundary, and that only a finite quantity of oil can be definitely trapped. Moreover, this quantity is  determined only by the capillary pressure curves and the difference between the volume mass of both phases, and does not depend on $u_0$. The model presented here, with total flow-rate $q$ equal to zero, do not allow this phenomenon, and all the oil present in $\O_1$ at the initial time remains trapped in $\O_1$ for all $t\ge 0$.

\bibliographystyle{plain}
\bibliography{ccances}

\end{document}